# Existence and Uniqueness of Solution for Linear Complementarity Problem in Contact Mechanics

**Jiamin Xu*, Nazli Demirer**, Vy Pho**, He Zhang***, Kaixiao Tian***, Ketan Bhaidasna**, Robert Darbe**, Dongmei Chen***

**The University of Texas at Austin, Austin, TX 78721 USA (e-mail: jiaminxu@utexas.edu, dmchen@me.utexas.edu).*

*** Halliburton, Houston, TX 77032 USA (e-mail: Nazli.Demirer@halliburton.com, vy.pho@halliburton.com, Robert.Darbe@halliburton.com, ketan.bhaidasna@halliburton.com )*

**** Halliburton, Singapore, SG 639940 Singapore (e-mail: He.Zhang@halliburton.com, Kaixiao.Tian@halliburton.com)*

**Abstract**: Although a unique solution is guaranteed in the Linear complementarity problem (LCP) when the matrix $\boldsymbol{M}$ is positive definite, practical applications often involve cases where $\boldsymbol{M}$ is only positive semi-definite, leading to multiple possible solutions. However, empirical observations suggest that uniqueness can still emerge under certain structural conditions on the vector $\boldsymbol{q}$. Motivated by an unresolved problem in nonlinear modeling for beam contact in directional drilling, this paper systematically investigates conditions under which a unique solution exists for LCPs with positive semi-definite matrices. We provide a rigorous proof demonstrating the existence and uniqueness of the solution for this specific case and extend our findings to establish a generalized framework applicable to broader classes of LCPs. This framework enhances the understanding of LCP uniqueness conditions and provides theoretical guarantees for solving real-world problems where positive semi-definite matrices arise.

## 1. INTRODUCTION

The Linear complementarity problem (LCP) involves determining a vector in a finite dimensional real vector space that meets specific inequality constraints. Formally, given a vector $\boldsymbol{q} \in \mathbb{R}^n$ and a matrix $\boldsymbol{M} \in \mathbb{R}^{n \times n}$, the task is to find a vector or a set of vectors $\boldsymbol{z} \in \mathbb{R}^n$ that satisfies the following conditions:

$$\boldsymbol{z} \geq 0 \tag{1.1}$$

$$\boldsymbol{q} + \boldsymbol{M}\boldsymbol{z} \geq 0 \tag{1.2}$$

$$\boldsymbol{z}^T(\boldsymbol{q} + \boldsymbol{M}\boldsymbol{z}) = 0 \tag{1.3}$$

Alternatively, the objective may be to prove that no such vector $\boldsymbol{z}$ exists. This LCP can be characterized by the pair $(\boldsymbol{q}, \boldsymbol{M})$.

**Remark 1**. In this paper, bold lowercase symbols represent vectors, bold uppercase symbols denote matrices, and unbold lowercase symbols indicate scalars. Moreover, the inequality for vectors means element wise. For example, a vector $\boldsymbol{z} \in \mathbb{R}^n$ with $\boldsymbol{z} \geq \boldsymbol{0}$ means every element in $\boldsymbol{z}$ is greater or equal to zero.

The Linear Complementarity Problem (LCP) has found broad applications in fields such as mathematics, operations research, computer science, game theory, economics, finance, and engineering. When the matrix $\boldsymbol{M}$ in the formulated LCP is positive definite, a unique solution exists [1], and efficient algorithms, such as the pivoting algorithm [2], can be employed to solve the problem. However, this is not always the case. In many practical applications, the matrix $\boldsymbol{M}$ is positive semi-definite [3–6], leading to multiple possible solutions, which can pose challenges in determining the appropriate solution.

In certain applications, however, even when $\boldsymbol{M}$ is positive semi-definite, numerous simulations suggest that a unique solution is highly likely due to the inherent structure of the vector $\boldsymbol{q}$ in equation (1.2) [7–10]. This observation motivates the theoretical investigation of the conditions under which the LCP, with $\boldsymbol{M}$ positive semi-definite, yields a unique solution when $\boldsymbol{q}$ exhibits specific structures.

Motivated by an unresolved problem [11,12] in nonlinear modeling for beam contact problem—where the problem is formulated as an LCP with a positive semi-definite $\boldsymbol{M}$—this paper systematically proves the existence and uniqueness of the solution for that specific case. By addressing this issue step by step, we not only resolve the original problem but also establish a more general framework for LCPs that ensures uniqueness when $\boldsymbol{M}$ is positive semi-definite and $\boldsymbol{q}$ has certain structures. Therefore, the primary contribution of this paper is the development of a generalized framework for ensuring the uniqueness of LCP solutions under certain conditions. This framework extends beyond contact modeling and can be applied to various fields where LCPs with positive semi-definite matrices arise.

The rest of the paper is organized as below: Section 2 provides a brief introduction to the unresolved problem in contact modeling for directional drilling, which also serves as an example demonstrating the application of the final theorem. Section 3 presents a step-by-step proof of the existence and uniqueness of the solution. Section 4 concludes the paper.

## 2. PROBLEM STATEMENT IN BEAM CONTACT

In the beam contact problem, an iterative algorithm [13,14] is commonly used to identify boundary conditions (i.e., contact points) before solving the system. Figure 1 provides a schematic representation of this process. This approach is adaptable for solving both linear and nonlinear beam contact problems, allowing contact to occur at any position along the beam. However, as the number of potential contact points increases, the computational burden becomes significant. Moreover, the possible presence of multiple solutions complicates the analysis, making it difficult to ascertain whether the obtained solution is correct and unique.

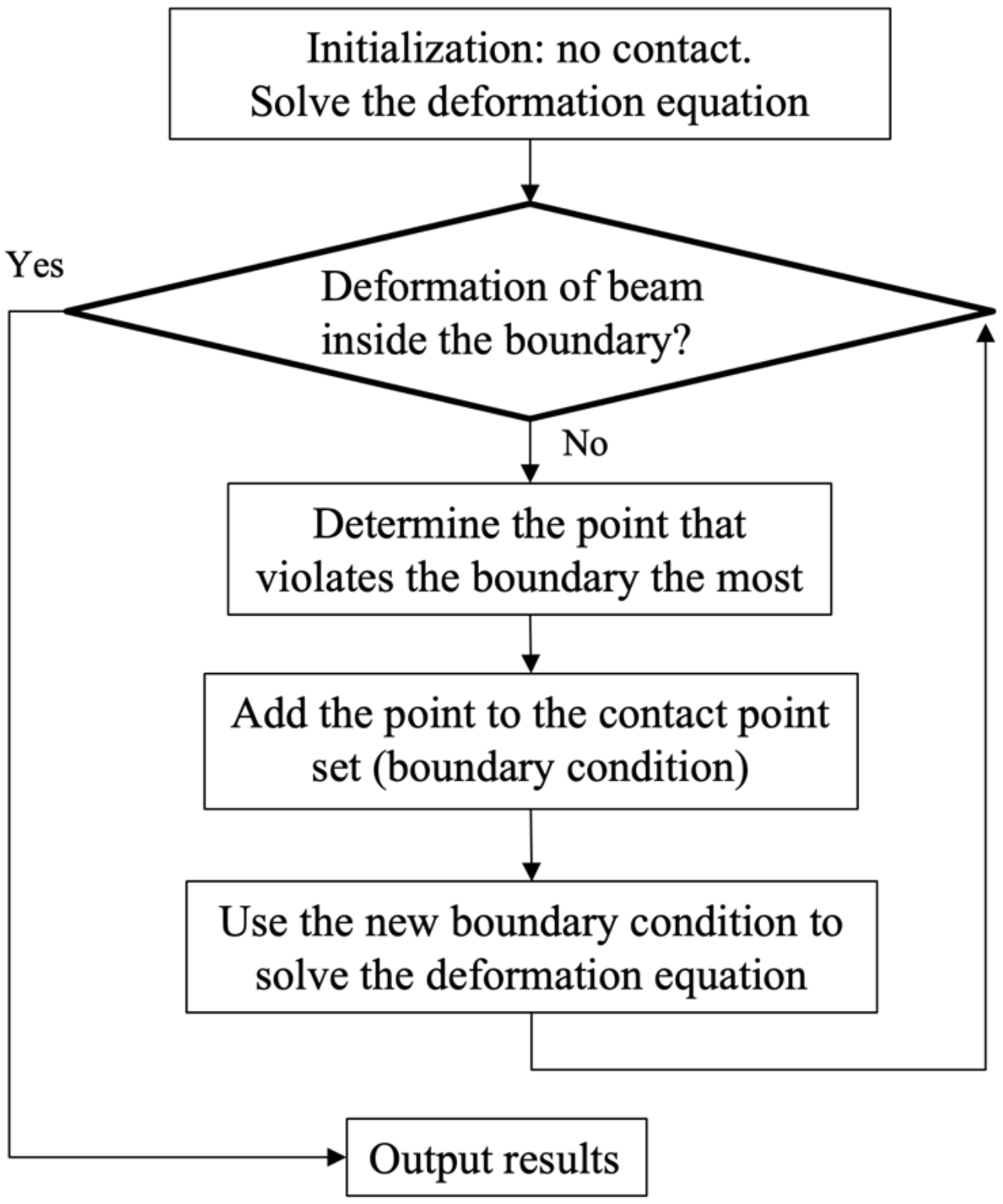

Figure 1. Schematics of the iteration algorithm for solving beam contact and deformation problem

For a 2D linear beam, when all potential contact positions are known in advance, the beam's deformation can be formulated as an LCP. This situation is common, particularly in applications like directional drilling [15–19], where the beam is intentionally designed with larger collars, known as stabilizers, to ensure that they make contact with the boundary rather than the adjacent main body of the beam. The schematics is shown in Figure 2 where all the stabilizers are the potential contact positions. Under different external inputs, such as the external lateral force, some stabilizers may be in contact with the boundary, while others may remain floating. The LCP is expressed as:

$$\underbrace{\begin{bmatrix}\boldsymbol{\gamma_l}\\ \boldsymbol{\gamma_u}\end{bmatrix}}_{\boldsymbol{w}} = \underbrace{\begin{bmatrix}\boldsymbol{K} & -\boldsymbol{K}\\ -\boldsymbol{K} & \boldsymbol{K}\end{bmatrix}}_{\boldsymbol{M}} \underbrace{\begin{bmatrix}\boldsymbol{F_l}\\ \boldsymbol{F_u}\end{bmatrix}}_{\boldsymbol{z}} + \underbrace{\begin{bmatrix}\tilde{\boldsymbol{q}} + \boldsymbol{y}^*\\ -\tilde{\boldsymbol{q}} + \boldsymbol{y}^*\end{bmatrix}}_{\boldsymbol{q}}, 0 \le \boldsymbol{w} \perp \boldsymbol{z} \ge 0 \quad (2.1)$$

$$\begin{cases} \boldsymbol{\gamma_l} = [\gamma_{l1}\ \gamma_{l2} \cdots \gamma_{ln}]^T \\ \boldsymbol{\gamma_u} = [\gamma_{u1}\ \gamma_{u2} \cdots \gamma_{un}]^T \\ \boldsymbol{F_l} = [F_{l1}\ F_{l2} \dots F_{ln}]^T \\ \boldsymbol{F_u} = [F_{u1}\ F_{u2} \dots F_{un}]^T \\ \tilde{\boldsymbol{q}} = [\tilde{q}_1\ \tilde{q}_2 \cdots \tilde{q}_n]^T \\ \boldsymbol{y}^* = [y_1^*\ y_2^* \dots y_n^*]^T \end{cases} \quad (2.2)$$

where $n$ is the number of possibly floating stabilizers; $\gamma_{li}$ ($\gamma_{ui}$) is the lower (upper) gap from the $ith$ stabilizer to the boundary; $F_{li}$ ($F_{ui}$) is the force when the $ith$ stabilizer is in contact with the lower (upper) boundary; $\tilde{q}_i$ is determined by the inputs to the beam such as the external axial and lateral forces; $y_i^* > 0$ represents the nominal gap of the $ith$ stabilizer when it is centered between the upper and lower boundaries; $\boldsymbol{K} \in \mathbb{R}^{n \times n}$ is the stiffness matrix that is symmetric and positive definite. Physically, this means if the lower (upper) contact force of a certain stabilizer, $F_{li}$ ($F_{ui}$) is greater than zero, such as the second (third) solid stabilizers from the right, then the corresponding lower (upper) gap $\gamma_{li}$ ($\gamma_{ui}$) must be zero. This applies to any stabilizers, leading to the $0 \le \boldsymbol{w} \perp \boldsymbol{z} \ge 0$ constraint in (2.1).

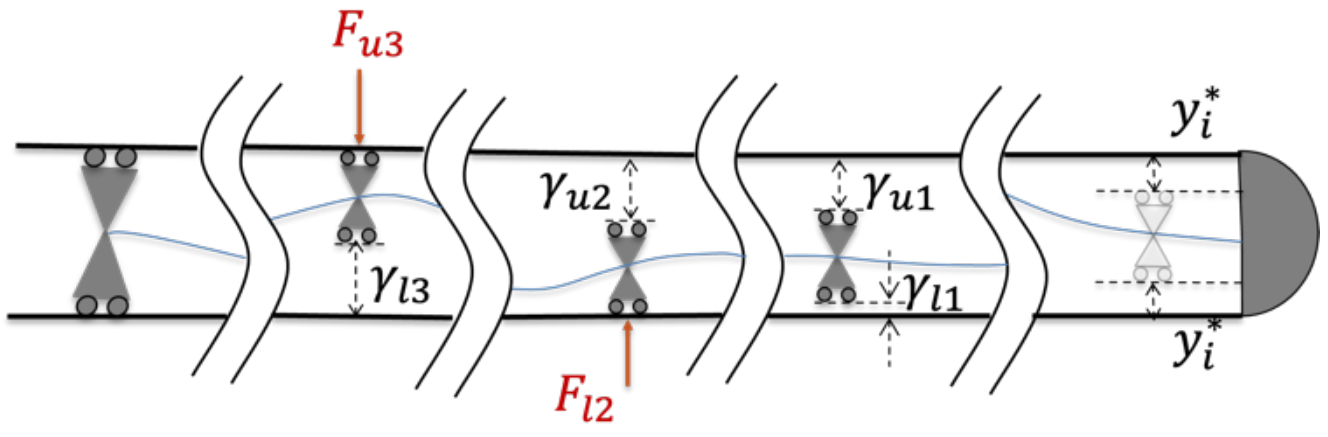

Figure 2. Schematics of the beam deformation and contact configuration in directional drilling.

In the original paper [11,12], this formulation became the core component of the state-of-the-art borehole propagation model in directional drilling, proposed by Shakib et al., which we will refer to as Shakib's model. Prior to this model [17–19], previous work assumed that all stabilizers were always in contact with the borehole, i.e., $\gamma_{li} = \gamma_{ui} = 0$ for any $i$. This assumption made the modeled beam excessively rigid and difficult to deform, resulting in an overall dynamic model that was less accurate compared to Shakib's model, which incorporates the formulation in Equation (2). Experimental validation has demonstrated the high accuracy of Shakib's model [10], leading to further research, including online estimation and control techniques [9]. However, one fundamental issue has remained unsolved for six years since the introduction of Shakib's model—the uniqueness of the solution, as noted in Shakib's paper. In the subsequent thesis [20], Shakib attempted to address this issue by using a variable transformation method to prove uniqueness for the special case of a single stabilizer (i.e., $n = 1$). However, the general case where $n > 1$ remains unresolved.

In the next section, we systematically prove the uniqueness of the solution without requiring any variable transformation in Equation (2). Additionally, we develop a generalized framework that ensures the uniqueness of the solution for broader applications.

## 3. PROOF OF EXISTENCE AND UNIQUENESS

### *3.1 Supporting definitions and theorems.*

*Definition 1.* A vector $\boldsymbol{z}$ that meets the conditions in inequalities (1.1) and (1.2) is termed feasible. We describe the LCP $(\boldsymbol{q}, \boldsymbol{M})$ as feasible when a feasible vector $\boldsymbol{z}$ exists.

*Remark 2.* According to Definition 1, if there exists a feasible $\boldsymbol{z}$ that satisfies equation (1.3), it represents a solution to the LCP $(\boldsymbol{q}, \boldsymbol{M})$. This condition is commonly referred to as 'solvable' in the literature; in this paper, however, we will simply state that a solution exists.

*Theorem 1* (Theorem 3.1.2 in [1]). Let $\boldsymbol{M}$ be a positive semi-definite matrix. If the LCP $(\boldsymbol{q}, \boldsymbol{M})$ is feasible, then there exists at least one solution.

*Theorem 2* (Theorem 3.1.7 in [1]). Let $\boldsymbol{M}$ be a symmetric and positive semi-definite matrix and $\boldsymbol{q} \in \mathbb{R}^n$ be arbitrary, if $\boldsymbol{z}^{(1)}$ and $\boldsymbol{z}^{(2)}$ are two solutions of $(\boldsymbol{q}, \boldsymbol{M})$, then the following hold:

1. $\left(\boldsymbol{z}^{(1)}\right)^{\mathrm{T}}\left(\boldsymbol{q} + \boldsymbol{M}\boldsymbol{z}^{(2)}\right) = \left(\boldsymbol{z}^{(2)}\right)^{\mathrm{T}}\left(\boldsymbol{q} + \boldsymbol{M}\boldsymbol{z}^{(1)}\right) = 0$
2. $\boldsymbol{M}\boldsymbol{z}^{(1)} = \boldsymbol{M}\boldsymbol{z}^{(2)}$

*Theorem 3* (Theorem 3.1.8 in [1]). Let $\boldsymbol{M} \in \mathbb{R}^{n\times n}$ and $\boldsymbol{q} \in \mathbb{R}^n$ be given, if for any two solutions $\boldsymbol{z}^{(1)}$ and $\boldsymbol{z}^{(2)}$ of $(\boldsymbol{q}, \boldsymbol{M})$, equation $\left(\boldsymbol{z}^{(1)}\right)^{\mathrm{T}}\left(\boldsymbol{q} + \boldsymbol{M}\boldsymbol{z}^{(2)}\right) = \left(\boldsymbol{z}^{(2)}\right)^{\mathrm{T}}\left(\boldsymbol{q} + \boldsymbol{M}\boldsymbol{z}^{(1)}\right) = 0$ holds, then the solution set of $(\boldsymbol{q}, \boldsymbol{M})$ is convex.

*Collaroy 1.* Combining *Theorem 2* and *Theorem 3*, it can be shown that if the LCP $(\boldsymbol{q}, \boldsymbol{M})$ is feasible and $\boldsymbol{M}$ is positive semi-definite, the solution set of the LCP $(\boldsymbol{q}, \boldsymbol{M})$ is convex.

### *3.2 Existence of solution.*

In this section, we used *Theorem 1* to prove that an LCP with the form of (2) has at least one solution. The two pre-requisites of *Theorem 1*, positive semi-definiteness of $\boldsymbol{M}$ and the feasibility of LCP $(\boldsymbol{M}, \boldsymbol{q})$, are given by *Theorem 4* and *Theorem 5*, respectively. The corresponding proof are given here as well.

*Theorem 4.* Let $\boldsymbol{K} \in \mathbb{R}^{n\times n}$ be a symmetric positive definite matrix. If a matrix $\boldsymbol{M}$ has the form of:

$$\boldsymbol{M} = \begin{bmatrix} \boldsymbol{K} & -\boldsymbol{K} \\ -\boldsymbol{K} & \boldsymbol{K} \end{bmatrix}$$

Then $\boldsymbol{M}$ is positive semi-definite.

*Proof:*

Let $\boldsymbol{x} = [\boldsymbol{u}\ \boldsymbol{v}]^T$, then

$$\boldsymbol{x}^{\mathrm{T}}\boldsymbol{M}\boldsymbol{x} = \boldsymbol{u}^{\mathrm{T}}\boldsymbol{K}\boldsymbol{u} - 2\boldsymbol{u}^{\mathrm{T}}\boldsymbol{K}\boldsymbol{v} + \boldsymbol{v}^{\mathrm{T}}\boldsymbol{K}\boldsymbol{v} = (\boldsymbol{u} - \boldsymbol{v})^{\mathrm{T}}\boldsymbol{K}(\boldsymbol{u} - \boldsymbol{v})$$

Given $K$ is positive definite, we have

$$\begin{Bmatrix} \boldsymbol{u} \neq \boldsymbol{v} & \boldsymbol{x}^{\mathrm{T}}\boldsymbol{M}\boldsymbol{x} > 0 \\ \boldsymbol{u} = \boldsymbol{v} & \boldsymbol{x}^{\mathrm{T}}\boldsymbol{M}\boldsymbol{x} = 0 \end{Bmatrix} \Rightarrow \boldsymbol{x}^{\mathrm{T}}\boldsymbol{M}\boldsymbol{x} \geq 0 \text{ for any } \boldsymbol{x}$$

$\Rightarrow \boldsymbol{M}$ is positive semi-definite. □

**Theorem 5**. An LCP with the form of

$$\underbrace{\begin{bmatrix} \boldsymbol{\gamma}_l \\ \boldsymbol{\gamma}_u \end{bmatrix}}_{\boldsymbol{w}} = \underbrace{\begin{bmatrix} \boldsymbol{K} & -\boldsymbol{K} \\ -\boldsymbol{K} & \boldsymbol{K} \end{bmatrix}}_{\boldsymbol{M}} \underbrace{\begin{bmatrix} \boldsymbol{F}_l \\ \boldsymbol{F}_u \end{bmatrix}}_{\boldsymbol{z}} + \underbrace{\begin{bmatrix} \tilde{\boldsymbol{q}} + \boldsymbol{y}^* \\ -\tilde{\boldsymbol{q}} + \boldsymbol{y}^* \end{bmatrix}}_{\boldsymbol{q}}$$

where $\boldsymbol{\gamma}_l, \boldsymbol{\gamma}_u, \boldsymbol{F}_l, \boldsymbol{F}_u, \tilde{\boldsymbol{q}} \in \mathbb{R}^{\boldsymbol{n}}$, $\boldsymbol{y}^* \in \mathbb{R}^{\boldsymbol{n}} > 0$ and $\boldsymbol{K} \in \mathbb{R}^{n\times n}$ is nonsingular, has at least one solution.

*Proof:*

Given $\boldsymbol{K}$ is nonsingular, we can always find two vectors $\boldsymbol{F}_l \geq \boldsymbol{0}$ and $\boldsymbol{F}_u \geq \boldsymbol{0}$ that satisfy

$$\boldsymbol{F}_l - \boldsymbol{F}_u = -\boldsymbol{K}^{-1}(\tilde{\boldsymbol{q}} + \boldsymbol{y}^*)$$

Subsequently we have:

$$\begin{Bmatrix} \boldsymbol{z} \geq \boldsymbol{0} \\ \begin{Bmatrix} \boldsymbol{\gamma}_l = \boldsymbol{0} \\ \boldsymbol{\gamma}_l = 2\boldsymbol{y}^* > \boldsymbol{0} \end{Bmatrix} \Rightarrow \boldsymbol{w} \geq \boldsymbol{0} \end{Bmatrix} \Rightarrow \text{The above LCP } (\boldsymbol{M}, \boldsymbol{q}) \text{ is}$$

feasible according to *Definition 1.*

Since $\boldsymbol{K} \in \mathbb{R}^{n\times n}$ is symmetric and positive definite, it is nonsingular, which will make the LCP $(\boldsymbol{M}, \boldsymbol{q})$ in (2) be feasible as just proved above. Then, according to *Theorem 4*, the matrix $M$ is positive semi-definite. Combining the abovementioned two conditions and using *Theorem 1* and *Collaroy 1*, we conclude that the LCP $(\boldsymbol{M}, \boldsymbol{q})$ in (2) has at least one solution and the solution set is convex. □

We will prove the solution is unique in the following section, i.e., the convex solution set is just one point.

### *3.3 Uniqueness of solution.*

In this section, we proved the uniqueness of the solution by exploring the specific structure of the LCP $(\boldsymbol{M}, \boldsymbol{q})$ in (2). Firstly, for all the solutions of (2), we have $\boldsymbol{F}_l^T\boldsymbol{F}_u = 0$ always be true as illustrated by *Theorem 6*. Then, we utilized *Theorem 2.2* and *Collaroy 1* to prove the uniqueness of the solution, summarized by *Theorem 7.*

*Theorem 6.* For all the solutions, $\boldsymbol{w}$ and $\boldsymbol{z}$ of (2), the equality $\boldsymbol{F}_l^T\boldsymbol{F}_u = 0$ always holds.

*Proof:*

Suppose $w$ and $z$ is one of the solutions of (2), then we have

$$\begin{cases} \boldsymbol{F}_l \geq \boldsymbol{0}, & (3.1) \\ \boldsymbol{F}_u \geq \boldsymbol{0}, & (3.2) \\ \boldsymbol{\gamma}_l = \boldsymbol{K}(\boldsymbol{F}_l - \boldsymbol{F}_u) + \tilde{\boldsymbol{q}} + \boldsymbol{y}^* \geq \boldsymbol{0}, & (3.3) \\ \boldsymbol{\gamma}_u = -\boldsymbol{K}(\boldsymbol{F}_l - \boldsymbol{F}_u) - \tilde{\boldsymbol{q}} + \boldsymbol{y}^* = 2\boldsymbol{y}^* - \boldsymbol{\gamma}_l \geq \boldsymbol{0}, & (3.4) \\ \boldsymbol{w}^T\boldsymbol{z} = 0 & (3.5) \end{cases}$$

From now on we use contradiction for the proof. Suppose $\boldsymbol{F}_l^T\boldsymbol{F}_u \neq 0$, then given

$\begin{Bmatrix}(3.1)\\(3.2)\end{Bmatrix} \Rightarrow$ there exists at least one $i = 1,2, \dots n$ such that $\begin{Bmatrix}F_{li} > 0\\F_{ui} > 0\end{Bmatrix}$, where $i$ is the index of the element in a vector of dimension $n$.

Subsequently, given (3.5), we have

$$\begin{Bmatrix}F_{l_i}\gamma_{l_i} = 0\\F_{u_i}\gamma_{u_i} = 0\end{Bmatrix} \Rightarrow \begin{Bmatrix}\gamma_{l_i} = 0\\\gamma_{u_i} = 0\end{Bmatrix} \Rightarrow \gamma_{l_i} + \gamma_{u_i} = 0\ .$$

However, according to (3.3) and (3.4), the following equality holds:

$$\boldsymbol{\gamma_l} + \boldsymbol{\gamma_u} = \boldsymbol{2y}^* \Rightarrow \gamma_{l_i} + \gamma_{u_i} = 2y_i^* > 0$$

This contradicts with the above equality $\gamma_{l_i} + \gamma_{u_i} = 0$ which is a consequence of $\boldsymbol{F_l^T F_u} \neq \boldsymbol{0}$. Thus, $\boldsymbol{F_l^T F_u} = 0$. □

*Theorem 6* can be easily understood through Figure 2. For example, if the lower force at one stabilizer $F_{li}$ is greater than zero, then the upper gap $\gamma_{li}$ has to be greater than zero, leading to no upper force there ($F_{li} = 0$). This is true for any stabilizers, thus $\boldsymbol{F_l^T F_u} = 0$.

*Theorem 7.* The LCP in *Theorem 5* also always has a unique solution.

**Proof**.

Suppose the above LCP $(\boldsymbol{M}, \boldsymbol{q})$ has at least two distinct solutions $\boldsymbol{z}^{(1)}$ and $\boldsymbol{z}^{(2)}$, then according to *Theorem 2.2* we have the below equality relationship

$$\boldsymbol{z}^{(2)} = \boldsymbol{z}^{(1)} + \mathcal{N}(\boldsymbol{M})\boldsymbol{x}$$

where $\mathcal{N}(\boldsymbol{M})$ is the null space of $\boldsymbol{M}$ and $\boldsymbol{x} \in \mathbb{R}^n$ is an arbitrary vector with $\|\boldsymbol{x}\|_2^2 \neq 0$. Given $\boldsymbol{K} \in \mathbb{R}^{n\times n}$ is symmetric positive definite and the structure of $\boldsymbol{M}$ in *Theorem 5*, the null space of $\boldsymbol{M}$ can be expressed as

$$\mathcal{N}(\boldsymbol{M}) = \begin{bmatrix}\boldsymbol{I}\\\boldsymbol{I}\end{bmatrix}$$

where $\boldsymbol{I} \in \mathbb{R}^{n\times n}$ is an identity matrix. Therefore, we have

$$\boldsymbol{z}^{(2)} = \begin{bmatrix}\boldsymbol{F_l}^{(2)}\\\boldsymbol{F_u}^{(2)}\end{bmatrix} = \boldsymbol{z}^{(1)} + \begin{bmatrix}\boldsymbol{x}\\\boldsymbol{x}\end{bmatrix} = \begin{bmatrix}\boldsymbol{F_l}^{(1)}\\\boldsymbol{F_u}^{(1)}\end{bmatrix} + \begin{bmatrix}\boldsymbol{x}\\\boldsymbol{x}\end{bmatrix} \tag{4}$$

Based on *Collaroy 1*, any vector $\boldsymbol{z}^{(3)}$ that is a convex combination of $\boldsymbol{z}^{(1)}$ and $\boldsymbol{z}^{(2)}$ is also a solution of the LCP $(\boldsymbol{M}, \boldsymbol{q})$:

$$\begin{aligned}\boldsymbol{z}^{(3)} &= \begin{bmatrix}\boldsymbol{F_l}^{(3)}\\\boldsymbol{F_u}^{(3)}\end{bmatrix}\\ &= \lambda\boldsymbol{z}^{(1)} + (1-\lambda)\boldsymbol{z}^{(2)}\\ &= \lambda\boldsymbol{z}^{(1)} + (1-\lambda)(\boldsymbol{z}^{(1)} + \begin{bmatrix}\boldsymbol{x}\\\boldsymbol{x}\end{bmatrix})\\ &= \boldsymbol{z}^{(1)} + (1-\lambda)\begin{bmatrix}\boldsymbol{x}\\\boldsymbol{x}\end{bmatrix}\\ &= \begin{bmatrix}\boldsymbol{F_l}^{(1)}\\\boldsymbol{F_u}^{(1)}\end{bmatrix} + (1-\lambda)\begin{bmatrix}\boldsymbol{x}\\\boldsymbol{x}\end{bmatrix}\end{aligned} \tag{5}$$

where $\lambda \in (0,1)$. According to *Theorem 6*, $\boldsymbol{z}^{(k)}, k = 1,2,3$, all needs to satisfy $\boldsymbol{F_l}^{(k)T}\boldsymbol{F_u}^{(k)} = 0$. Substitute this expression for (4) into the equality we have

$$\boldsymbol{F_l}^{(2)T}\boldsymbol{F_u}^{(2)} = \underbrace{\boldsymbol{F_l}^{(1)T}\boldsymbol{F_u}^{(1)}}_{0} + \boldsymbol{x}^T\left(\boldsymbol{F_l}^{(1)} + \boldsymbol{F_u}^{(1)}\right) + \underbrace{\boldsymbol{x}^T\boldsymbol{x}}_{\|\boldsymbol{x}\|_2^2} = 0$$

$$\Rightarrow\ \boldsymbol{x}^T\left(\boldsymbol{F_l}^{(1)} + \boldsymbol{F_u}^{(1)}\right) = -\|\boldsymbol{x}\|_2^2 \tag{6}$$

By doing the same operation for $\boldsymbol{z}^{(3)}$ we have

$$\begin{aligned}&\boldsymbol{F_l}^{(3)T}\boldsymbol{F_u}^{(3)}\\ &= \underbrace{\boldsymbol{F_l}^{(1)T}\boldsymbol{F_u}^{(1)}}_{0} + (1-\lambda)\boldsymbol{x}^T\left(\boldsymbol{F_l}^{(1)} + \boldsymbol{F_u}^{(1)}\right) + (1-\lambda)^2\underbrace{\boldsymbol{x}^T\boldsymbol{x}}_{\|\boldsymbol{x}\|_2^2}\\ &= (1-\lambda)\boldsymbol{x}^T\left(\boldsymbol{F_l}^{(1)} + \boldsymbol{F_u}^{(1)}\right) + (1-\lambda)^2\|\boldsymbol{x}\|_2^2 = 0\end{aligned} \tag{7}$$

Substituting (6) into (7):

$$\begin{aligned}&-(1-\lambda)\|\boldsymbol{x}\|_2^2 + (1-\lambda)^2\|\boldsymbol{x}\|_2^2\\ &= \underbrace{\lambda}_{\neq 0}\underbrace{(\lambda - 1)}_{\neq 0}\|\boldsymbol{x}\|_2^2 = 0\\ &\Rightarrow\ \|\boldsymbol{x}\|_2^2 = 0\end{aligned}$$

This contradicts with the requirement $\|\boldsymbol{x}\|_2^2 \neq 0$ for having two distinct solutions.

Therefore, the LCP $(\boldsymbol{M}, \boldsymbol{q})$ has a unique solution. Furthermore, given that the existence of solution proved in Section 3.2, we can conclude that the LCP $(\boldsymbol{M}, \boldsymbol{q})$ always has and only has one solution. □

*2.4 Generalization.*

*Theorem 8.* An LCP $(\boldsymbol{M}, \boldsymbol{q})$ that has the following structure always has a unique solution $\boldsymbol{z} = [\boldsymbol{z}_1^T\ \boldsymbol{z}_2^T \dots\ \boldsymbol{z}_t^T]^T$:

$$\boldsymbol{M} = \begin{bmatrix}\boldsymbol{M}_{11} & & & \\ \boldsymbol{M}_{21} & \boldsymbol{M}_{22} & & \\ \vdots & \vdots & \ddots & \\ \boldsymbol{M}_{t1} & \boldsymbol{M}_{t2} & \cdots & \boldsymbol{M}_{tt}\end{bmatrix}$$

$$\boldsymbol{q} = \begin{bmatrix}\boldsymbol{q}_1\\\boldsymbol{q}_2\\\vdots\\\boldsymbol{q}_t\end{bmatrix}$$

With

$$\left\{\begin{aligned}&\boldsymbol{M}_{ii} = \begin{bmatrix}\boldsymbol{K}_i & -\boldsymbol{K}_i\\-\boldsymbol{K}_i & \boldsymbol{K}_i\end{bmatrix}\\ &\boldsymbol{M}_{ij} = \begin{bmatrix}\widetilde{\boldsymbol{K}}_{ij} & -\widetilde{\boldsymbol{K}}_{ij}\\-\widetilde{\boldsymbol{K}}_{ij} & \widetilde{\boldsymbol{K}}_{ij}\end{bmatrix}\\ &\boldsymbol{q}_i = \begin{bmatrix}\boldsymbol{q}_{i1}\\\boldsymbol{q}_{i2}\end{bmatrix}\\ &\boldsymbol{q}_{i1} + \boldsymbol{q}_{i2} = \boldsymbol{q}_{ic} > 0\end{aligned}\right\}\ i,j = 1,2,\dots,t, and\ j < i$$

where $\boldsymbol{K}_i \in \mathbb{R}^{n_i\times n_i}$ are symmetric and positive definite matrices; $\widetilde{\boldsymbol{K}}_{ij} \in \mathbb{R}^{n_i\times n_{i-1}}$ are arbitrary matrices (does not have to be square); $\boldsymbol{q}_{i1}, \boldsymbol{q}_{i2}, \boldsymbol{q}_{ic} \in \mathbb{R}^n$ and $\boldsymbol{z}_i \in \mathbb{R}^{2n_i}$.

*Proof:*

Starting from the solution $\boldsymbol{z}_1$ in $\boldsymbol{z}$, which corresponds to $\boldsymbol{M}_{11}$ and is decoupled from the remaining. We have:

$$\underbrace{\begin{bmatrix}\boldsymbol{w}_{11}\\ \boldsymbol{w}_{12}\end{bmatrix}}_{\boldsymbol{w}_1} = \underbrace{\begin{bmatrix}\boldsymbol{K}_1 & -\boldsymbol{K}_1\\ -\boldsymbol{K}_1 & \boldsymbol{K}_1\end{bmatrix}}_{\boldsymbol{M}_{11}} \underbrace{\begin{bmatrix}\boldsymbol{z}_{11}\\ \boldsymbol{z}_{12}\end{bmatrix}}_{\boldsymbol{z}_1} + \underbrace{\begin{bmatrix}\boldsymbol{q}_{11}\\ \boldsymbol{q}_{12}\end{bmatrix}}_{\boldsymbol{q}_1} \quad (8)$$

which has the same form of (2) by noticing that (2) is equivalent to $\boldsymbol{q}_{11} = \widetilde{\boldsymbol{q}} + \boldsymbol{y}^*, \boldsymbol{q}_{12} = -\widetilde{\boldsymbol{q}} + \boldsymbol{y}^*, \boldsymbol{q}_{11} + \boldsymbol{q}_{12} = 2\boldsymbol{y}^* = \boldsymbol{q}_{1c} > 0$ in (8). Then according to *Theorem 7*, (8) always has a unique solution $\boldsymbol{z}_1$.

Moving to $\boldsymbol{z}_2$, we have

$$\underbrace{\begin{bmatrix}\boldsymbol{w}_{21}\\ \boldsymbol{w}_{22}\end{bmatrix}}_{\boldsymbol{w}_2} = \underbrace{\begin{bmatrix}\boldsymbol{K}_2 & -\boldsymbol{K}_2\\ -\boldsymbol{K}_2 & \boldsymbol{K}_2\end{bmatrix}}_{\boldsymbol{M}_{22}} \underbrace{\begin{bmatrix}\boldsymbol{z}_{21}\\ \boldsymbol{z}_{22}\end{bmatrix}}_{\boldsymbol{z}_2} + \underbrace{\begin{bmatrix}\boldsymbol{q}_{21}\\ \boldsymbol{q}_{22}\end{bmatrix}}_{\boldsymbol{q}_2} + \underbrace{\begin{bmatrix}\widetilde{\boldsymbol{K}}_{21} & -\widetilde{\boldsymbol{K}}_{21}\\ -\widetilde{\boldsymbol{K}}_{21} & \widetilde{\boldsymbol{K}}_{21}\end{bmatrix}}_{\boldsymbol{M}_{21}} \underbrace{\begin{bmatrix}\boldsymbol{z}_{11}\\ \boldsymbol{z}_{12}\end{bmatrix}}_{\boldsymbol{z}_1}$$

$$= \underbrace{\begin{bmatrix}\boldsymbol{K}_2 & -\boldsymbol{K}_2\\ -\boldsymbol{K}_2 & \boldsymbol{K}_2\end{bmatrix}}_{\boldsymbol{M}_{22}} \underbrace{\begin{bmatrix}\boldsymbol{z}_{21}\\ \boldsymbol{z}_{22}\end{bmatrix}}_{\boldsymbol{z}_2} + \underbrace{\begin{bmatrix}\underbrace{\boldsymbol{q}_{21} + \widetilde{\boldsymbol{K}}_{21}\boldsymbol{z}_{11} - \widetilde{\boldsymbol{K}}_{21}\boldsymbol{z}_{12}}_{\widehat{\boldsymbol{q}}_{21}}\\ \underbrace{\boldsymbol{q}_{22} - \widetilde{\boldsymbol{K}}_{21}\boldsymbol{z}_{11} + \widetilde{\boldsymbol{K}}_{21}\boldsymbol{z}_{12}}_{\widehat{\boldsymbol{q}}_{22}}\end{bmatrix}}_{\widehat{\boldsymbol{q}}_2} \quad (9)$$

where $\widehat{\boldsymbol{q}}_2$ follows

$$\widehat{\boldsymbol{q}}_{21} + \widehat{\boldsymbol{q}}_{22} = \boldsymbol{q}_{21} + \boldsymbol{q}_{22} = \boldsymbol{q}_{2c}$$

Then (9) also always has a unique solution $\boldsymbol{z}_2$. The same operation can be done for all the other $z_i$ and we can conclude the LCP $(\boldsymbol{M}, \boldsymbol{q})$ in *Theorem 8* always has a unique solution. ☐

## 4. CONCLUSIONS

This paper systematically addresses the problem of solution uniqueness in LCPs with positive semi-definite matrices. While classical results guarantee uniqueness only for positive definite matrices, our work demonstrates that uniqueness can still be ensured under specific structural conditions on the vector $\boldsymbol{q}$. By resolving an open problem in beam contact modeling, we establish a generalized theoretical framework that extends to other LCP formulations with positive semi-definite matrices.

A key contribution of this work is Theorem 8, which provides a formal guarantee for the uniqueness of LCP solutions under particular conditions. This theorem has broad applications beyond the beam contact problem and can be largely utilized in contact [21] and grasp [5] modeling in robotics, where complementarity conditions often arise in the analysis of force closure and object manipulation. Additionally, it can be applied in mechanical systems, control theory, and optimization, offering a rigorous foundation for addressing uniqueness issues in these domains.

## ACKNOWLEDGEMENT

This research was supported by Halliburton Energy Services. The authors would like to thank Halliburton for the continuous support and the permission to publish this paper.